                    \def\version{1 March, 2011}                                %

\documentclass[reqno,11pt]{amsart} 
\usepackage{amsmath} 
\usepackage[mathscr]{eucal}
\usepackage{amssymb} 
\usepackage{srcltx} 
\numberwithin{equation}{section}
 
\newenvironment{proofsect}[1] 
{\vskip0.1cm\noindent{\bf #1.}\hskip0.5cm}

\newfam\Bbbfam 
\font\tenBbb=msbm10 
\font\sevenBbb=msbm7 
\font\fiveBbb=msbm5 
\textfont\Bbbfam=\tenBbb 
\scriptfont\Bbbfam=\sevenBbb 
\scriptscriptfont\Bbbfam=\fiveBbb

\newcommand{\R}     {\mathbb{R}} 
\newcommand{\Z}     {\mathbb{Z}} 
\newcommand{\N}     {\mathbb{N}} 
\renewcommand{\P}   {\mathbb{P}} 
 
\newcommand{\E}     {\mathbb{E}} 
 
\newcommand{\smfrac}[2]{\textstyle{\frac {#1}{#2}}}

\def\1{{\mathchoice {1\mskip-4mu\mathrm l}      
{1\mskip-4mu\mathrm l} 
{1\mskip-4.5mu\mathrm l} {1\mskip-5mu\mathrm l}}} 

\newcommand{\ssup}[1] {{{\scriptscriptstyle{({#1}})}}} 
\def\comment#1{} 
\newtheoremstyle{thm}{2ex}{2ex}{\itshape\rmfamily}{} 
{\bfseries\rmfamily}{}{1.7ex}{} 
 
\newtheoremstyle{rem}{1.3ex}{1.3ex}{\rmfamily}{} 
{\itshape\rmfamily}{}{1.5ex}{}

\newtheorem{theorem}{Theorem}[section] 
\newtheorem{lemma}[theorem]{Lemma}

\theoremstyle{definition}

\renewcommand{\d}{{\rm d}} 
 
\newcommand{\eps}{\varepsilon}

\newcommand{\e}   {{\operatorname e }}
\newcommand{\g}   {{\operatorname {maxexc}}}  
 
\begin{document} 
 
\title[The longest excursion of a random interacting polymer]{\large The longest excursion\\ of a random interacting polymer} 
 
\author[Janine K\"ocher and Wolfgang 
        K{\"o}nig]{} 
\maketitle

\thispagestyle{empty} 
 
\centerline {\sc Janine K\"ocher\footnote{Institute for Mathematics, TU Berlin, Str.~des 17.~Juni 136, 10623 Berlin, Germany, {\tt JanineKoecher@aol.com} and {\tt koenig@math.tu-berlin.de}} and Wolfgang
K{\"o}nig\footnotemark[1]$^,$\footnote{Weierstrass Institute Berlin, Mohrenstr.~39, 10117 Berlin, {\tt koenig@wias-berlin.de}}}
\vspace{0.4cm}
\renewcommand{\thefootnote}{}

\centerline{\small(\version)} 
\vspace{.5cm} 

\begin{quote} 
{\small {\bf Abstract:}} We consider a random $N$-step polymer under the influence of an attractive interaction with the origin and derive a limit law -- after suitable shifting and norming -- for the length of the longest excursion towards the Gumbel distribution. The embodied law of large numbers in particular implies that the longest excursion is of order $\log N$ long. The main tools are taken from extreme value theory and renewal theory.
\end{quote}



\setcounter{section}{0} 
\section{Introduction and main results}\label{Intro}

\noindent Let $(S_n)_{n\in\N_0}$ be a random walk on the lattice $\Z^d$ starting at the origin and having steps of mean zero. By $\P$ and $\E$ we denote the corresponding probability and expectation, respectively. We conceive the walk $(n,S_n)_{n=0,\dots,N}$ as an $N$-step polymer in the $(d+1)$-dimensional space. We introduce an attractive interaction with the origin by introducing the Gibbs measure $\P_{\beta,N}$ via the density
\begin{equation}
\frac{\d\P_{\beta,N}}{\d \P}=\frac{\e^{\beta L_N}} {Z_{\beta,N}}\qquad \mbox{with }Z_{\beta,N}=\E\big[ \e^{\beta L_N}\big],
\end{equation}
where $\beta\in(0,\infty)$ is a parameter and 
\begin{equation}
L_N=|\{k\in\{1,\dots,N\}\colon S_k=0\}|
\end{equation}
denotes the walker's local time at the origin, i.e., the number of returns to the origin. The properties of the polymer under $\P_{\beta,N}$ have been studied a lot \cite{dH09,G07}. In particular, the free energy
\begin{equation}
 F(\beta)=\lim_{N\to\infty}\frac1N\log Z_{\beta,N}\in(0,\beta)
\end{equation}
has been shown to exist and to be positive and strictly increasing in $\beta$. Furthermore, it has been shown that the polymer is localised in the sense that $L_N$ is of order $N$ under $\P_{\beta, N}$, and the density of the set of hits of the origin has been characterised. In particular, the constrained version, i.e., the polymer under 
\begin{equation}
\P_{\beta, N}^{\ssup{\rm c}}(\cdot)=\frac1{Z_{\beta,N}^{\ssup{\rm c}}}\E\big[\e^{\beta L_N}\1\{\,\cdot\,\}\1_\{S_N=0\}\big],\qquad \mbox{where }Z_{\beta,N}^{\ssup{\rm c}}=\E\big[\e^{\beta L_N}\1_\{S_N=0\}\big],
\end{equation}
has been studied.

In this paper, we consider the length of the longest excursion of the polymer under $\P_{\beta, N}^{\ssup{\rm c}}$. To introduce this object, we denote by $\tau=\{\tau_i\colon i\in\N_0\}$ the set of return times to the origin, where
\begin{equation}
\tau_0=0\qquad\mbox{and, inductively, }\tau_{i+1}=\inf\{n>\tau_i\colon S_n=0\},\quad i\in\N_0.
\end{equation}
Then $\P_{\beta, N}^{\ssup{\rm c}}$ is the conditional distribution of the polymer given $\{N\in\tau\}$. The length of the longest excursion is now given as
\begin{equation}
\g_N=\max\{\tau_{i}-\tau_{i-1}\colon i\in\N,\tau_{i}\leq N\}.
\end{equation}
According to \cite[Theorem 7.3]{dH09}, $\g_N$ is of order $\log N$ under $\P_{\beta, N}^{\ssup{\rm c}}$, in the sense that the distribution of $\g_N/\log N$ under $\P_{\beta, N}^{\ssup{\rm c}}$ is tight in $N$. The proof gives the upper bound $2/F(\beta)$, which is not sharp, as we will see below. It is the main goal of this note to derive not only the law of large numbers for $\g_N$, but also a non-trivial limit law for $\g_N$ after suitable shifting, in the spirit of extreme value theory.

To formulate our main result, we need to fix our assumptions first. 

\noindent{\bf Assumption ($\boldsymbol\tau$).} There are $D\in(0,\infty)$ and $\alpha\in(1,\infty)$ such that
$$
K(n):=\P(\tau_1=n)\sim D n^{-\alpha},\qquad n\to\infty.
$$

This assumption is fulfilled for most of the aperiodic random walks $(S_n)_{n\in\N_0}$ under consideration in the literature. For random walks with period $p\in\N$, one has to work with $K(pn)$ instead of $K(n)$ and with $pN$-step polymers and obtains analogous results. Assumption ($\tau$) can be relaxed with the help of slowly varying functions, on cost of a more cumbersome formulation and proof of the main result.

The main result of this paper is the following.

\begin{theorem}\label{thm-Main} Suppose that Assumption ($\tau$) is satisfied, and fix $\beta\in(0,\infty)$. Then, as $N\to\infty$, the distribution of 
\begin{equation}
 F(\beta)\g_N-\log\frac N{\mu_\beta}+\alpha \log\log\frac N{\mu_\beta}-C
\end{equation}
under $\P_{\beta, N}^{\ssup{\rm c}}$ weakly converges towards the standard Gumbel distribution, where
\begin{equation}
\mu_\beta=\e^\beta\sum_{n\in\N}n K(n) \e^{-n F(\beta)}\qquad \mbox{and}\qquad C=\log\Big(F(\beta)^\alpha D\frac{\e^{\beta-F(\beta)}}{1-\e^{-F(\beta)}}\Big).
\end{equation}
\end{theorem}

Explicitly, it is stated that, for any $x\in \R$,
\begin{equation}
\lim_{N\to\infty}\P_{\beta, N}^{\ssup{\rm c}}\big(\g_N\leq \gamma_x(N/\mu_\beta)\big)=\e^{-\e^{-x}},\qquad\mbox{where }\gamma_x(N)=\frac {x+C+\log N-\alpha \log\log N}{F(\beta)}.
\end{equation}
In particular, we have the law of large numbers: $\g_N/\log N\to 1/F(\beta)$ in $\P_{\beta,N}^{\ssup{\rm c}}$-probability as $N\to\infty$.

\section{The proof}

\noindent It is well-known that the free energy $F(\beta)$ is characterised by the equation
\begin{equation}
\e^\beta=\sum_{n\in\N}K(n)\e^{-n F(\beta)},
\end{equation}
and that it actually holds that $Z^{\ssup{\rm c}}_{\beta,N}\sim \e^{N F(\beta)}\frac1{\mu_\beta}$ as $N\to\infty$. In particular, $F(\beta)$ is also the exponential rate of $Z^{\ssup{\rm c}}_{\beta,N}$. The first step, which is basic to all investigations of the polymer, is a change of measure to the measure $Q_\beta$, under which the excursion lengths $T_k=\tau_{k+1}-\tau_k$, are i.i.d.~in $k\in\N_0$ with distribution
$$
Q_\beta(T_1=n)=\e^{-\beta}K(n)\e^{-n F(\beta)},\qquad n\in\N.
$$
Since $\g_N$ is measurable with respect to the family of the $T_k$'s, it is easy to see from the technique explained in \cite[p.~9]{G07} that
\begin{equation}
\P_{\beta,N}^{\ssup{\rm c}}(\g_N\leq \gamma_N)\sim \mu_\beta Q_\beta(\g_N\leq \gamma_N,N\in\tau),\qquad N\to\infty,
\end{equation}
for any choice of the sequence $(\gamma_N)_{N\in\N}$, where $\mu_\beta=\sum_{n\in\N}n Q_\beta(T_1=n)\in[1,\infty)$ is the expectation of the length of the first excursion under $Q_\beta$. Introducing 
\begin{equation}\label{MNsigmaNdef}
M_n=\max_{k=1}^n T_k\qquad \mbox{and}\qquad \sigma_N=\inf\{k\in\N\colon \tau_k\geq N\},
\end{equation}
we see that $\g_N= M_{\sigma_N}$ on $\{N\in\tau\}$ for any $N\in\N$. (Note that $\sigma_N=L_N$ on the event $\{N\in\tau\}$.) Hence, Theorem~\ref{thm-Main} is equivalent to 
\begin{equation}\label{sufficient}
\lim_{N\to\infty}Q_\beta(M_{\sigma_N}\leq \gamma_x(N/\mu_\beta), N\in\tau )=\frac 1{\mu_\beta}\e^{-\e^{-x}},\qquad x\in\R.
\end{equation}
The proof of this consist of a combination of three fundamental ingredients:
\begin{enumerate}
\item[(1)] an {\it extreme value theorem} for $M_n$ under $Q_\beta$,

\item[(2)] a {\it law of large numbers} for $\sigma_N$ under $Q_\beta$,

\item[(3)] a {\it renewal theorem} for $\tau$ under $Q_\beta$.
\end{enumerate}

Items (2) and (3) are immediate: We have from renewal theory that $\sigma_N/N\to 1/\mu_\beta$ in $Q_\beta$-probability and $\lim_{N\to\infty}Q_\beta(N\in\tau)=1/\mu_\beta$. The first item needs a bit more care:

\begin{lemma}\label{lem-extreme}
$$
\lim_{N\to\infty}Q_\beta(M_N\leq \gamma_x(N))=\e^{-\e^{-x}},\qquad x\in\R.
$$
\end{lemma}

\begin{proofsect}{Proof} Note that $M_N$ is the maximum of $N$ independent random variables with the same distribution as $T_1=\tau_1$ under $Q_\beta$. Observe that the tails of this distribution are given by
$$
\begin{aligned}
Q_\beta(\tau_1 > k)&=\e^{\beta}\sum_{n> k}K(n)\e^{-n F(\beta)}\sim \e^\beta D \sum_{n> k}n^{-\alpha}\e^{-n F(\beta)}\\
&=\e^\beta D \e^{-k F(\beta)}k^{-\alpha}\sum_{n\in\N}(1+\smfrac nk)^{-\alpha}\e^{-n F(\beta)}\\
&\sim \e^{-k F(\beta)}k^{-\alpha}D\frac {\e^{\beta-F(\beta)}}{1-\e^{-F(\beta)}},\qquad k\to\infty,
\end{aligned}
$$
where in the last step we used the monotonous convergence theorem and the geometric series. Hence, replacing $k$ by $\gamma_x(N)$, we see that, as $N\to\infty$,
$$
\begin{aligned}
Q_\beta(\tau_1 > \gamma_x(N))&\sim \e^{-\gamma_x(N) F(\beta)}\gamma_x(N)^{-\alpha}D\frac {\e^{\beta-F(\beta)}}{1-\e^{-F(\beta)}}\\
&=\frac 1N \e^{-C-x}(\log N)^\alpha\Big(\frac{x+C+\log N-\alpha \log\log N}{F(\beta)}\Big)^{-\alpha}
\e^C F(\beta)^{-\alpha}\\
&\sim \frac{\e^{-x}}N .
\end{aligned}
$$
From this  the assertion easily follows.
\qed
\end{proofsect}

Hence, Theorem~\ref{thm-Main} is easily seen to follow from the above three ingredients, as soon as one shows that $\sigma_N$ may asymptotically be replaced by $N/\mu_\beta$ and that the two events in \eqref{sufficient} are asymptotically independent. This is what we show now. First we show that $M_{\sigma_N}$ and $M_{N/\mu_\beta}$ have the same limiting distribution.

\begin{lemma}\label{cor-MsigmaN}
$$
\lim_{N\to\infty}Q_\beta(M_{\sigma_N}\leq \gamma_x(N/\mu_\beta))=\e^{-\e^{-x}},\qquad x\in\R.
$$
\end{lemma}

\begin{proofsect}{Proof}
The upper bound is proved as follows. Fix a small $\eps>0$, then we have, as $N\to\infty$,
\begin{equation}\label{Corproof}
\begin{aligned}
Q_\beta(M_{\sigma_N}\leq \gamma_x(N/\mu_\beta))
&\leq Q_\beta\Big(M_{\sigma_N}\leq \gamma_x(N/\mu_\beta),\sigma_N\geq \frac N{\mu_\beta+\eps}\Big)+Q_\beta\Big(\sigma_N<\frac N{\mu_\beta+\eps}\Big)\\
&\leq Q_\beta\big(M_{N/(\mu_\beta+\eps)}\leq \gamma_x(N/\mu_\beta)\big)+o(1).
\end{aligned}
\end{equation}
Observe that, as $N\to\infty$,
$$
\begin{aligned}
\gamma_x(N/\mu_\beta)-\gamma_x(N/(\mu_\beta+\eps))
&=\frac 1{F(\beta)}\log(1+{\smfrac\eps{\mu_\beta}})+\frac\alpha{F(\beta)}\log\frac{\log N-\log(\mu_\beta+\eps)}{\log N-\log\mu_\beta}\\
&=\frac 1{F(\beta)}\log(1+{\smfrac\eps{\mu_\beta}})+o(1).
\end{aligned}
$$
Hence, we may replace, as an upper bound, $\gamma_x(N/\mu_\beta)$ on the right of \eqref{Corproof} by $\gamma_{x+B\eps}(N/(\mu_\beta+\eps))$ for some suitable $B\in\R$, use Lemma~\ref{lem-extreme} for $N$ replaced by $N/(\mu_\beta+\eps)$ and $x$ replaced by $x+B\eps$ and make $\eps\downarrow 0$ in the end. This shows that the upper bound of the assertion holds. The lower bound is proved in the same way.
\qed
\end{proofsect}

\begin{proofsect}{Proof of Theorem~\ref{thm-Main}}
It is convenient to introduce a Markov chain $(Y_n)_{n\in\N_0}$ with 
$$
Y_n=\big(Y^{\ssup 1}_n,Y^{\ssup 2}_n\big)=\big(T_{\sigma_n},\tau_{\sigma_n}-n\big)
$$ 
on the state space $I=\{(i,j)\in\N\times\N_0\colon j\leq i\}$, where we recall \eqref{MNsigmaNdef}. In words, the first component is the size of the step over $n$, and the last is the size of the overshoot. This Markov chain is ergodic and positiv recurrent with invariant distribution $\pi(i,j)=Q_\beta(\tau_1=i)/\mu_\beta$ for $(i,j)\in I$. We denote by $\widetilde Q_{i,j}$ the distribution of this chain given that it starts in $Y_0=(i,j)$; note that $Q_\beta=\widetilde Q_{i,0}$ with an unspecified value of $i$, which we put equal to 1 by default. The event $\{N\in\tau\}$ is identical to $\{Y_N^{\ssup 2}=0\}=\{Y_N\in\N\times\{0\}\}$; by ergodicity, its probability under $\widetilde Q_{i,j}$ converges, as $N\to\infty$, to $\pi(\N\times \{0\})=\frac1{\mu_\beta}$, for any $(i,j)\in I$, which is one way to prove the renewal theorem.

Now let $\eps>0$ be given. Pick $K_\eps\in\N$ so large that $\pi(I_{K_\eps}^{\rm c})<\eps/2$, where $I_k=\{(i,j)\in I\colon i\leq k\}$ for any $k\in\N$. Furthermore, pick $R_\eps\in\N$ with $R_\eps>K_\eps$ so large that $\widetilde Q_{i,j}(R_\eps\in \tau)\leq \frac 1{\mu_\beta}+\eps$ for any $(i,j)\in I_{K_\eps}$. Now pick $N_\eps\in\N$ so large that $N_\eps>R_\eps$ and $\widetilde Q_{1,0}(Y_{N-R_\eps}\in I_{K_\eps}^{\rm c})<\pi(I_{K_\eps}^{\rm c})+\eps/2$ for any $N\geq N_\eps$. The latter is possible, since $\widetilde Q_{1,0}(Y_{N-R_\eps}\in I_{K_\eps}^{\rm c})=1-\widetilde Q_{1,0}(Y_{N-R_\eps}\in I_{K_\eps})$ converges towards $1-\pi(I_{K_\eps})=\pi(I_{K_\eps}^{\rm c})$ as $N\to\infty$ by ergodicity.

Recall that we only have to prove \eqref{sufficient}. We calculate, with the help of the Markov property at time $N-R_\eps$, for $N>N_\eps$,
$$
\begin{aligned}
Q_\beta&(M_{\sigma_N}\leq \gamma_x(N/\mu_\beta), N\in\tau )
=\widetilde Q_{1,0}\Big(\max_{k=1}^NY_k^{\ssup 1}\leq \gamma_x(N/\mu_\beta), Y_N^{\ssup 2}=0\Big)\\
&\leq \widetilde Q_{1,0}\Big(\max_{k=1}^{N-R_\eps} Y_k^{\ssup 1}\leq \gamma_x(N/\mu_\beta), Y_{N-R_\eps}\in I_{K_\eps},Y_N^{\ssup 2}=0\Big)+\widetilde Q_{1,0}(Y_{N-R_\eps}\in I_{K_\eps}^{\rm c})\\
&\leq \sum_{(i,j)\in I_{K_\eps}}\widetilde Q_{1,0}\Big(\max_{k=1}^{N-R_\eps}Y_k^{\ssup 1}\leq \gamma_x(N/\mu_\beta), Y_{N-R_\eps}=(i,j)\Big)\,\widetilde Q_{i,j}(Y_{R_\eps}^{\ssup 2}=0)+\pi(I_{K_\eps}^{\rm c})+\eps/2\\
&\leq \widetilde Q_{1,0}\Big(\max_{k=1}^{N-R_\eps}Y_k^{\ssup 1}\leq \gamma_x(N/\mu_\beta)\Big)(\smfrac 1{\mu_\beta}+\eps)+\eps\\
&\leq Q_\beta\big(M_{\sigma_{N-R_\eps}}\leq \gamma_x(N/\mu_\beta)\big)(\smfrac 1{\mu_\beta}+\eps)+\eps.
\end{aligned}
$$
Now apply Lemma~\ref{cor-MsigmaN} for $N$ replaced by $N-R_\eps$ and observe that $\lim_{N\to\infty}(\gamma_x(N/\mu_\beta)-\gamma_x((N-R_\eps)/\mu_\beta))=0$. Afterwards letting $\eps\downarrow0$ shows that the upper bound in \eqref{sufficient} holds. The proof of the corresponding lower bound is similar, and we omit it.
\qed
\end{proofsect}

\end{document}